\documentclass{amsart}
\usepackage{amsmath,amscd,amssymb}
\usepackage[all]{xy}

\allowdisplaybreaks
\numberwithin{equation}{section}
\theoremstyle{plain}
\newtheorem{thm}{Theorem}[section]
\newtheorem{cor}[thm]{Corollary}
\newtheorem{lem}[thm]{Lemma}
\newtheorem{prop}[thm]{Proposition}
\theoremstyle{definition}

\newtheorem{rem}[thm]{Remark}

\def\ee{\mathaccent23e}
\newcommand{\proj}{\operatorname{pr}}

\begin{document}

\title[Cellular decomposition and L-S cat of $Spin(7)$]
{On the cellular decomposition and the Lusternik-Schnirelmann
  category  of $Spin(7)$}
\author{Norio Iwase, Mamoru Mimura and Tetsu Nishimoto}
\address{Faculty of Mathematics,
  Kyushu University,
  Ropponmatsu Fukuoka 810-8560, Japan}
\email{iwase@math.kyushu-u.ac.jp}
\address{Department of Mathematics,
  Faculty of Science,
  Okayama University,
  3-1 Tsushima-naka, Okayama 700-8530, Japan}
\email{mimura@math.okayama-u.ac.jp}
\address{Department of Welfare Business,
  Kinki Welfare University, Fukusaki-cho,
  Hyogo 679-2217, Japan}
\email{nishimoto@kinwu.ac.jp}

\begin{abstract}
  We give a cellular decomposition of the compact connected Lie group
  $Spin(7)$.
  We also determine the L-S categories of $Spin(7)$ and $Spin(8)$.
\end{abstract}

\maketitle

\section{Introduction}

In this paper, we assume that a space has the homotopy type of a
CW-complex.
The Lusternik-Schnirelmann category ${\rm cat}\,X$ of a space $X$ is
the least integer $n$ such that $X$ is the union of $(n+1)$ open
subsets, each of which is contractible in $X$.
G.~Whitehead \cite{white1} showed that ${\rm cat}\,X \leq n$ if and only if the
diagonal map $\Delta_{n+1} : X \rightarrow \prod^{n+1} X$ is homotopic
to some composition map
\begin{equation*}
  X \longrightarrow T^{n+1}(X) \longrightarrow \prod^{n+1} X,
\end{equation*}
where $T^{n+1}(X)$ is the fat wedge and $T^{n+1}(X) \rightarrow
\prod^{n+1} X$ is the inclusion map.

The weak Lusternik-Schnirelmann category $w{\rm cat}\,X$ is the least
integer $n$ such that the reduced diagonal map $\bar{\Delta}_{n+1} : X
\rightarrow \wedge^{n+1} X$ is trivial.
Then it is easy to see that $w{\rm cat}\,X \leq {\rm cat}\,X$,
since $\wedge^{n+1} X = \prod^{n+1} X/T^{n+1}(X)$.

The strong Lusternik-Schnirelmann category ${\rm Cat}\,X$ is
the least integer $n$ such that there exist a space $X'$ which is
homotopy equivalent to $X$ and is covered by $(n+1)$ open subsets
contractible in themselves.
${\rm Cat}\,X$ is closely related with ${\rm cat}\,X$, and
Ganea and Takens \cite{takens} showed that
\begin{equation*}
  {\rm cat}\,X \leq {\rm Cat}\,X \leq {\rm cat}\,X + 1.
\end{equation*}
Ganea \cite{ganea} showed that ${\rm Cat}\,X$ is equal to the invariant
which is the least integer $n$ such that there is a cofibre sequence
\begin{equation*}
  A_i \longrightarrow X_{i-1} \longrightarrow X_i
\end{equation*}
where $X_0$ is a point and $X_n$ is homotopy equivalent to $X$.

The Lusternik-Schnirelmann category for some Lie groups are
determined, such as ${\rm cat}(U(n)) = n$ and ${\rm cat}(SU(n)) = n-1$
by Singhof \cite{s1}, ${\rm cat}(Sp(2)) = 3$ by Schweitzer \cite{sch},
${\rm cat}(Sp(3)) = 5$ by
Fern\'andez-Su\'arez, G\'omez-Tato, Strom and Tanr\'e \cite{fgst},
and Iwase and Mimura \cite{im},
${\rm cat}(SO(2)) = 1$, ${\rm cat}(SO(3)) = 3$, ${\rm cat}(SO(4)) =
4$, ${\rm cat}(SO(5)) = 8$
by James and Singhof \cite{js}.
Some general argument about the Lusternik-Schnirel-mann category
implies that ${\rm cat}(G_2) = 4$ (see for example \cite{im}).

As is well-known, we have the following isomorphisms:
\begin{equation*}
    Spin(3) \cong S^3, \quad  Spin(4) \cong S^3 \times S^3, \quad
    Spin(5) \cong Sp(2), \quad Spin(6) \cong SU(4).
\end{equation*}
Thus $Spin(7)$ is the first non-trivial case in determining the
cellular decomposition and the Lusternik-Schnirelmann category as well;
it is our purpose in this paper.

\begin{thm}
  \label{thm:main}
  We have $w{\rm cat}(Spin(7)) = {\rm cat}(Spin(7)) =
  {\rm Cat}(Spin(7)) = 5$.
\end{thm}

Since $Spin(8)$ is homeomorphic to $Spin(7) \times S^7$, we obtain the
following corollary.

\begin{cor}
  We have $w{\rm cat}(Spin(8)) = {\rm cat}(Spin(8)) =
  {\rm Cat}(Spin(8)) = 6$.
\end{cor}

The paper is organized as follows.
In Section 2 we give a cellular decomposition of $Spin(7)$ such that
$Spin(7)$ contains a subgroup $SU(4)$, which turns out to be useful
for determining the Lusternik-Schnirelmann category of $Spin(7)$.
In Section 3 we give a cone-decomposition of $SU(4)$, which gives rise 
to the Lusternik-Schnirelmann category of $Spin(7)$ in Section 4.

\section{The cellular decomposition of $Spin(7)$}

In this section, we use the notation in \cite{mn}.
Let $\mathfrak C$ be the Cayley algebra.
$SO(8)$ acts on $\mathfrak C$ naturally since $\mathfrak C \cong
\mathbb R^8$ as $\mathbb R$-module.
We regard $SO(7)$ as the subgroup of $SO(8)$ fixing $e_0$, the 
unit of $\mathfrak C$.
As is well known, the exceptional Lie group $G_2$ is defined by
\begin{equation*}
  G_2 = \{ g\in SO(7) \; | \; g(x)g(y) = g(xy), x,y \in \mathfrak C
  \} = {\rm Aut}(\mathfrak C).
\end{equation*}
According to \cite{yok3}, the group $Spin(7)$ is the set of the
elements $\tilde{g} \in SO(8)$ such that $g(x)\tilde{g}(y) =
\tilde{g}(xy)$ for any $x, y \in \mathfrak C$, where $g \in SO(7)$ is
uniquely determined by $\tilde{g}$:
\begin{equation*}
  Spin(7) = \{ \tilde{g} \in SO(8) \; | \; g(x)\tilde{g}(y)
  = \tilde{g}(xy), g \in SO(7), x,y \in \mathfrak C \}.
\end{equation*}
It is easy to see that $G_2$ is the subgroup of $Spin(7)$.
Observe that the algebra generated by $e_1$ in $\mathfrak C$ is
isomorphic to $\mathbb C$.
$SU(4)$ acts on $\mathfrak C$ naturally, since as $\mathbb C$-module
$\mathfrak C \cong \mathbb C^4$ whose basis is $\{ e_0, e_2, e_4,
e_6 \}$.
We regard $SU(3)$ as the subgroup of $SU(4)$ fixing $e_0$ and also as
the subgroup of $G_2$ fixing $e_1$.

Let $D^i$ be the $i$-dimensional disc.
We define four maps:
\begin{align*}
  & A : D^3 = \{ (x_1,x_2,x_3) \in \mathbb R^3 \; | \; x_1^2+x_2^2+x_3^2
  \leq 1 \}
  \longrightarrow SO(8), \\
  & B : D^2 = \{ (y_1,y_2) \in \mathbb R^2 \; | \; y_1^2+y_2^2 \leq 1 \}
  \longrightarrow SO(8), \\
  & C : D^1 = \{ z_1 \in \mathbb R \; | \; z_1^2 \leq 1 \}
  \longrightarrow SO(8), \\
  & D : D^2 = \{ (w_1,w_2) \in \mathbb R^2 \; | \;
  w_1^2+w_2^2 \leq 1 \} \longrightarrow SO(8), \\
\end{align*}
as follows:
\begin{align*}
  & A(x_1, x_2, x_3) =
  \begin{pmatrix}
    1 & & & & & & & \\
    & 1 & & & & & & \\
    & & 1 & & & & & \\
    & & & 1 & & & & \\
    & & & & 1-2X^2
    & - 2x_1X
    & -2x_2X
    & -2x_3X \\
    & & & & 2x_1X
    & 1-2X^2
    & 2x_3X
    & -2x_2X \\
    & & & & 2x_2X
    & -2x_3X
    & 1-2X^2
    & 2x_1X \\
    & & & & 2x_3X
    & 2x_2X
    & -2x_1X
    & 1-2X^2 \\
  \end{pmatrix}, \\
  & B(y_1, y_2) =
  \begin{pmatrix}
    1 & & & & & & & \\
    & 1 & & & & & & \\
    & & y_1 & -y_2 & -Y & 0 & & \\
    & & y_2 & y_1 & 0 & -Y & & \\
    & & Y & 0 & y_1 & y_2 & & \\
    & & 0 & Y & -y_2 & y_1 & & \\
    & & & & & & 1 & \\
    & & & & & & & 1 \\
  \end{pmatrix}, \\
  & C(z_1) =
  \begin{pmatrix}
    1 & & & & & & & \\
    & z_1 & 0 & -Z & & & & \\
    & 0 & 1 & 0 & & & & \\
    & Z & 0 & z_1 & & & & \\
    & & & & 1 & & & \\
    & & & & & z_1 & 0 & -Z \\
    & & & & & 0 & 1 & 0 \\
    & & & & & Z & 0 & z_1 \\
  \end{pmatrix}, \\
  & D(w_1, w_2) =
  \begin{pmatrix}
    w_1 & -w_2 & -W   &      &   &   &   & \\
    w_2 &  w_1 &      & -W   &   &   &   & \\
    W   &      &  w_1 &  w_2 &   &   &   & \\
        &  W   & -w_2 &  w_1 &   &   &   & \\
        &      &      &      & 1 &   &   & \\
        &      &      &      &   & 1 &   & \\
        &      &      &      &   &   & 1 & \\
        &      &      &      &   &   &   & 1 \\
  \end{pmatrix},
\end{align*}
where we put for simplicity
\begin{equation*}
  \begin{array}{ll}
  X=\sqrt{1-x_1^2-x_2^2-x_3^2}, & Y=\sqrt{1-y_1^2-y_2^2}, \\
  Z=\sqrt{1-z_1^2}, & W=\sqrt{1-w_1^2-w_2^2}.
  \end{array}
\end{equation*}

\begin{lem}
  \label{lem:in_spin7}
  The elements $A(x_1,x_2,x_3)$, $B(y_1,y_2)$, $C(z_1)$ and
  $D(w_1,w_2)$ belong to $Spin(7)$.
\end{lem}
\begin{proof}
  Apparently the elements $A(x_1,x_2,x_3)$, $B(y_1,y_2)$ and $C(z_1)$
  belong to $G_2$.
  In the proof, we denote $D(w_1,w_2)$ simply by $D$.
  Let $D'$ be the matrix
  \begin{equation*}
    \begin{pmatrix}
      1 &   &   &   &      &      &     & \\
        & 1 &   &   &      &      &     & \\
        &   & 1 &   &      &      &     & \\
        &   &   & 1 &      &      &     & \\
        &   &   &   &  w_1 & -w_2 & W   & \\
        &   &   &   &  w_2 &  w_1 &     & -W \\
        &   &   &   & -W   &      & w_1 & -w_2 \\
        &   &   &   &      &  W   & w_2 &  w_1 \\
    \end{pmatrix}.
  \end{equation*}
  Then we can show by a tedious calculation that $D'xDy=D(xy)$ for any
  $x, y \in \mathfrak C$, which gives us the result.
\end{proof}

Let $\varphi_3$, $\varphi_5$, $\varphi_6$ and $\varphi_7$ be maps
\begin{align*}
  & \varphi_3 : D^3 \longrightarrow Spin(7), \\
  & \varphi_5 : D^3 \times D^2 \longrightarrow Spin(7), \\
  & \varphi_6 : D^3 \times D^2 \times D^1 \longrightarrow Spin(7), \\
  & \varphi_7 : D^3 \times D^2 \times D^2 \longrightarrow Spin(7)
\end{align*}
respectively defined by the equalities
\begin{align*}
  & \varphi_3 ({\bf x}) = A({\bf x}), \\
  & \varphi_5 ({\bf x},{\bf y}) =
  B({\bf y})A({\bf x})B({\bf y})^{-1}, \\
  & \varphi_6 ({\bf x},{\bf y},{\bf z}) =
  C({\bf z})B({\bf y})A({\bf x})B({\bf y})^{-1}C({\bf z})^{-1}, \\
  & \varphi_7 ({\bf x},{\bf y},{\bf w}) =
  D({\bf w})B({\bf y})A({\bf x})B({\bf y})^{-1}D({\bf w})^{-1},
\end{align*}
where ${\bf x} = (x_1,x_2,x_3)$, ${\bf y} = (y_1,y_2)$, ${\bf z} =
(z_1)$ and ${\bf w} = (w_1,w_2)$.
We define sixteen cells $e^j$ for $j = 0, 3, 5, 6, 7, 8, 9, 10, 11, 12,
13, 14, 15, 16, 18, 21$ respectively as follows:
\begin{equation*}
  \begin{array}{lllll}
    e^0 = \{ 1 \}, & e^3 = {\rm Im} \; \varphi_3,
    & e^5 = {\rm Im} \; \varphi_5, & e^6 = {\rm Im} \; \varphi_6,
    & e^7 = {\rm Im} \; \varphi_7, \\
    e^8 = e^5e^3, & e^9 = e^6e^3, & e^{10} = e^7e^3,
    & e^{11} = e^6e^5, & e^{12} = e^7e^5, \\
    e^{13} = e^6e^7, & e^{14} = e^6e^5e^3, & e^{15} = e^7e^5e^3,
    & e^{16} = e^6e^7e^3, & e^{18} = e^6e^7e^5, \\
    e^{21} = e^6e^7e^5e^3.
  \end{array}
\end{equation*}

Let $S^7$ be the unit sphere of $\mathfrak C$.
Then we have a principal bundle over it:
\begin{equation*}
  SU(3) \longrightarrow SU(4) \overset{p_0}{\longrightarrow} S^7,
\end{equation*}
where $p_0(g) = ge_0$.

\begin{lem}
  Let $V^7 = D^3 \times D^2 \times D^2$.
  Then the composite map $p_0 \varphi_7 : (V^7, \partial V^7) \rightarrow
  (S^7, e_0)$ is a relative homeomorphism.
\end{lem}
\begin{proof}
  We express the map $(p_0 \varphi_7)|_{V^7\backslash \partial V^7}$
  as follows:
  \begin{equation*}
    \begin{pmatrix}
      a_0 \\
      a_1 \\
      a_2 \\
      a_3 \\
      a_4 \\
      a_5 \\
      a_6 \\
      a_7
    \end{pmatrix}
    = D({\bf w})B({\bf y})A({\bf x})B({\bf y})^{-1}D({\bf w})^{-1}e_0
    =
    \begin{pmatrix}
      1 - 2X^2Y^2W^2 \\
      2x_1XY^2W^2 \\
      2(w_1X-x_1w_2)XY^2W \\
      -2(w_2X+x_1w_1)XY^2W \\
      2(-y_1X+x_1y_2)XYW \\
      2(y_2X+x_1y_1)XYW \\
      2x_2XYW \\
      2x_3XYW
    \end{pmatrix}
  \end{equation*}
  and hence we have
  \begin{equation*}
    \begin{pmatrix}
      1-a_0 \\
      a_1 \\
      a_2 \\
      a_3 \\
      a_4 \\
      a_5 \\
      a_6 \\
      a_7
    \end{pmatrix}
    = 2XYW
    \begin{pmatrix}
      XYW \\
      x_1YW \\
      (w_1X-x_1w_2)Y \\
      -(w_2X+x_1w_1)Y \\
      -y_1X+x_1y_2 \\
      y_2X+x_1y_1 \\
      x_2 \\
      x_3
    \end{pmatrix}.
  \end{equation*}
  Since $X > 0$, $Y > 0$, $W > 0$ and $1-a_0 > 0$, an easy calculation
  as for the first component in the above equation gives the following
  equation:
  \begin{equation}
    \label{eq:subXYW}
    XYW = \frac{\sqrt{1-a_2}}{\sqrt2},
  \end{equation}
  from which we easily obtain
  \begin{equation}
    \label{eq:subx23}
    x_2 = \frac{a_6}{\sqrt{2(1-a_2)}}, \quad
    x_3 = \frac{a_7}{\sqrt{2(1-a_2)}}.
  \end{equation}
  Further we obtain three more equalities from the above equalities:
  \begin{align*}
    & (1-a_0)^2 + a_1^2 = 4X^2Y^4W^4(x_1^2 + X^2), \\
    & a_2^2 + a_3^2 = 4X^2Y^4W^2(w_1^2 + w_2^2)(x_1^2 + X^2) =
    4X^2Y^4W^2(1 - W^2)(x_1^2 + X^2), \\
    & a_4^2 + a_5^2 = 4X^2Y^2W^2(y_1^2 + y_2^2)(x_1^2 + X^2) =
    4X^2Y^2W^2(1 - Y^2)(x_1^2 + X^2).
  \end{align*}
  Using these three equalities, we obtain
  \begin{align}
    \label{eq:subY2}
    & Y^2 = \frac{(1-a_0)^2 + a_1^2 + a_2^2 + a_3^2}
    {(1-a_0)^2 + a_1^2 + a_2^2 + a_3^2 + a_4^2 + a_5^2}, \\
    \label{eq:subW2}
    & W^2 = \frac{(1-a_0)^2 + a_1^2}{(1-a_0)^2 + a_1^2 + a_2^2 + a_3^2}.
  \end{align}
  It follows from (\ref{eq:subXYW}), (\ref{eq:subY2}) and (\ref{eq:subW2})
  that
  \begin{equation}
    \label{eq:subX2}
    X^2 = \frac{(1-a_0)((1-a_0)^2 + a_1^2 + a_2^2 + a_3^2 + a_4^2 + a_5^2)}
    {2((1-a_0)^2 + a_1^2)}.
  \end{equation}
  It follows also from (\ref{eq:subx23}) and
  (\ref{eq:subX2}) that
  \begin{equation}
    \label{eq:subx12}
    x_1^2 = \frac{a_1^2((1-a_0)^2 + a_1^2 + a_2^2 + a_3^2 + a_4^2 + a_5^2)}
    {2(1-a_0)((1-a_0)^2 + a_1^2)}.
  \end{equation}
  Since $Y$, $W$, $X$ are positive, (\ref{eq:subY2}),
  (\ref{eq:subW2}), (\ref{eq:subX2}) imply respectively
  \begin{align}
    \label{eq:subY}
    & Y = \frac{\sqrt{(1-a_0)^2 + a_1^2 + a_2^2 + a_3^2}}
    {\sqrt{(1-a_0)^2 + a_1^2 + a_2^2 + a_3^2 + a_4^2 + a_5^2}}, \\
    \label{eq:subW}
    & W = \frac{\sqrt{(1-a_0)^2 + a_1^2}}
    {\sqrt{(1-a_0)^2 + a_1^2 + a_2^2 + a_3^2}}, \\
    \label{eq:subX}
    & X = \frac{\sqrt{(1-a_0)((1-a_0)^2 + a_1^2 + a_2^2 + a_3^2
    + a_4^2 + a_5^2)}}{\sqrt{2((1-a_0)^2 + a_1^2)}}.
  \end{align}
  Since the signs of $x_1$ and $a_1$ are the same,
  (\ref{eq:subx12}) implies that
  \begin{equation}
    \label{eq:subx1}
    x_1 = \frac{a_1\sqrt{(1-a_0)^2 + a_1^2 + a_2^2 + a_3^2 + a_4^2 + a_5^2}}
    {\sqrt{2(1-a_0)((1-a_0)^2 + a_1^2)}}.
  \end{equation}
  Now we determine $y_1$;
  we have
  \begin{equation*}
    -a_4X + a_5x_1 = 2XYW(x_1^2+X^2)y_2.
  \end{equation*}
  Substituting the equations (\ref{eq:subXYW}), (\ref{eq:subX}) and
  (\ref{eq:subx1}) in the above equation, we obtain
  \begin{equation}
    \label{eq:suby1}
    y_1 = \frac{a_1a_5-(1-a_0)a_4}
    {\sqrt{((1-a_0)^2+a_1^2)((1-a_0)^2+a_1^2+a_2^2+a_3^2+a_4^2+a_5^2)}}.
  \end{equation}
  We determine $y_2$;
  we have
  \begin{equation*}
    a_4x_1 + a_5X = 2XYW(x_1^2+X^2)y_2.
  \end{equation*}
  Substituting the equations (\ref{eq:subXYW}), (\ref{eq:subX}) and
  (\ref{eq:subx1}) in the above equation, we obtain
  \begin{equation}
    \label{eq:suby2}
    y_2 = \frac{a_1a_4+(1-a_0)a_5}
    {\sqrt{((1-a_0)^2+a_1^2)((1-a_0)^2+a_1^2+a_2^2+a_3^2+a_4^2+a_5^2)}}.
  \end{equation}
  We determine $w_1$;
  we have
  \begin{equation*}
    a_2X - a_3x_1 = 2XY^2W(x_1^2+X^2)w_1.
  \end{equation*}
  Substituting the equations (\ref{eq:subXYW}), (\ref{eq:subY}),
  (\ref{eq:subX}) and (\ref{eq:subx1}) in the above equation, we obtain
  \begin{equation}
    \label{eq:subw1}
    w_1 = \frac{(1-a_0)a_2-a_1a_3}
    {\sqrt{((1-a_0)^2+a_1^2)((1-a_0)^2+a_1^2+a_2^2+a_3^2)}}.
  \end{equation}
  Finally we determine $w_2$;
  we have
  \begin{equation*}
    - a_2x_1 - a_3X = 2XY^2W(x_1^2+X^2)w_2.
  \end{equation*}
  Substituting the equations (\ref{eq:subXYW}), (\ref{eq:subY}),
  (\ref{eq:subX}) and (\ref{eq:subx1}) in the above equation, we obtain
  \begin{equation}
    \label{eq:subw2}
    w_2 = \frac{-a_1a_2-(1-a_0)a_3}
    {\sqrt{((1-a_0)^2+a_1^2)((1-a_0)^2+a_1^2+a_2^2+a_3^2)}}.
  \end{equation}
  Thus we have expressed $x_1, x_2, x_3, y_1, y_2, w_1, w_2$ in terms
  of $a_0, \cdots, a_7$, that is, the inverse map has been constructed,
  which completes the proof.
\end{proof}

In a similar way to that of Section 3 of \cite{mn}, we can obtain the
following theorem, which is essentially the same as Yokota's
decomposition \cite{yok1}.

\begin{prop}
  \label{prop:su4}
  $e^0 \cup e^3 \cup e^5 \cup e^7 \cup e^8 \cup e^{10} \cup e^{12} \cup
  e^{15}$ thus obtained is a cellular decomposition of $SU(4)$.
\end{prop}
\begin{proof}
  First we show that $\ee^i \cap \ee^j
  = \emptyset$ if $i \not= j$.
  We consider the following three cases:

  (1) For the case where $i, j \in \{0, 3, 5, 8 \}$;
  both cells $e^i$ and $e^j$ are in $SU(3)$ and $e^0 \cup e^3 \cup
  e^5 \cup e^8$ is a cellular decomposition of $SU(3)$.
  Then we have $\ee^i \cap \ee^j
  = \emptyset$ if $i \not= j$.

  (2) For the case where $i \in \{0, 3, 5, 8 \}$ and $j \in \{7, 10,
  12, 15 \}$;
  we have $p_0(\ee^i) = \{e_0\}$ and $p_0(\ee^j)
  = S^7 \backslash \{e_0\}$.
  Then we have $\ee^i \cap \ee^j
  = \emptyset$.

  (3) For the case where $i, j \in \{7, 10, 12, 15 \}$;
  suppose that $A \in \ee^i \cap \ee^j$.
  Since $\ee^i = \ee^7\ee^{i-7}$ and $\ee^j = \ee^7\ee^{j-7}$, we
  can put $A = A_1A_2 = A'_1A'_2$ where $A_1, A'_1 \in \ee^7$, $A_2
  \in \ee^{i-7}$ and $A'_2 \in \ee^{j-7}$.
  We have $A_1 = A'_1$, since $p_0(A_1) = p_0(A_1A_2) = p_0(A'_1A'_2)
  = p_0(A'_1)$ and $p_0|_{\ee^7}$ is monic.
  Then we have $A_2 = A'_2$ and the first case shows that $i-7 =
  j-7$, that is, $i=j$.
  Thus $\ee^i \cap \ee^j = \emptyset$ if $i \not=j$.

  Next, we will check that the boundaries of the cells are included in
  the lower dimensional cells.
  In the proof of Proposition 3.2 \cite{mn}, it is proved that
  the boundaries $\dot e^3$, $\dot e^5$ and $\dot e^8$ are included in
  the lower dimensional cells.
  Observe that the boundary $\dot e^7$ is the union of the following
  three sets:
  \begin{align*}
    & \{ DBAB^{-1}D^{-1} \; | \; A \in A(\dot D^3), B \in B(D^2),
    D \in D(D^2) \}, \\
    & \{ DBAB^{-1}D^{-1} \; | \; A \in A(D^3), B \in B(\dot D^2),
    D \in D(D^2) \}, \\
    & \{ DBAB^{-1}D^{-1} \; | \; A \in A(D^3), B \in B(D^2),
    D \in D(\dot D^2) \}.
  \end{align*}
  The first set contains only the identity element, since $A$ is the
  identity element.
  It is easy to see that the second set is contained in $e^3$ and
  that the third set is contained in $e^5$.
  We have $\dot e^{10} = e^7\dot e^3 \cup \dot e^7e^3 \subset e^7e^0 \cup
  e^5e^3 = e^7 \cup e^8$.
  We also have $\dot e^{12} = \dot e^7e^5 \cup e^7\dot e^5 \subset
  e^5e^5 \cup e^7e^3 = e^8 \cup e^{10}$, and
  $\dot e^{15} = \dot e^7e^5e^3 \cup e^7\dot e^5e^3 \cup
  e^7e^5\dot e^3 \subset e^5e^5e^3 \cup e^7e^3e^3 \cup e^7e^5 = e^8
  \cup e^{10} \cup e^{12}$.

  Finally, we will show that the inclusion map $e^0 \cup e^3 \cup e^5
  \cup e^7 \cup e^8 \cup e^{10} \cup e^{12} \cup e^{15} \rightarrow SU(4)$
  is epic.
  Let $g \in SU(4)$.
  If $p_0(g) = e_0$, then $g$ is contained in $SU(3) = e^0 \cup e^3
  \cup e^5 \cup e^8$.
  Suppose that $p_0(g) \not= e_0$.
  There is an element $h\in e^7$ such
  that $p_0(h) = p_0(g)$.
  Thus we have $h^{-1}g \in SU(3) = e^0 \cup e^3 \cup e^5 \cup e^8$,
  since $p_0(h^{-1}g) = e_0$.
  Therefore we have $g \in h(e^0 \cup e^3 \cup e^5 \cup e^8) \subset
  e^0 \cup e^3 \cup e^5 \cup e^7 \cup e^8 \cup e^{10} \cup e^{12} \cup
  e^{15}$.
\end{proof}

\begin{rem}
  \label{rem:e}
  (1) We regard $SO(6)$ as the subgroup of $SO(7)$ fixing $e_1$.
  Let $\pi : Spin(6) \rightarrow SO(6)$ be the double covering.
  Then, according to the Proof of Lemma \ref{lem:in_spin7}, $\pi(SU(4))
  \subset SO(6)$ so that $\pi|_{SU(4)} : SU(4) \rightarrow SO(6)$ is the
  double covering.

  (2) For $1 \leq n \leq 3$, the subcomplex
  $e^0 \cup e^3 \cup \cdots \cup e^{2n+1}$ is
  homeomorphic to $\Sigma \mathbb CP^n$, which consists of the elements
  \begin{equation*}
    A
    \begin{pmatrix}
      1 & & & \\
      & 1 & & \\
      & & 1 & \\
      & & & e^{2i\theta}
    \end{pmatrix}
    A^{-1}
    \begin{pmatrix}
      1 & & & \\
      & 1 & & \\
      & & 1 & \\
      & & & e^{-2i\theta}
    \end{pmatrix}
  \end{equation*}
  for any elements $A$ in $SU(n+1)$.
  Moreover, according to Proposition 2.6 of Chapter IV of \cite{st}, we
  have $e^{2i+1}e^{2j+1} \subset e^{2j+1}e^{2i+1}$ for $i < j$; in
  fact we have $e^{2i+1}e^{2j+1} = e^{2j+1}e^{2i+1}$ (see \cite{yok3}).
\end{rem}

Let $S^6$ be the unit sphere of $\mathbb R^7$ whose basis $\{ e_i \; |
\; 1 \leq i \leq 7 \}$.
We consider the following diagram
\begin{equation*}
  \begin{CD}
    SU(3) @>>> G_2 @>>> S^6 \\
    @VVV @VVV @| \\
    SU(4) @>>> Spin(7) @>p>> S^6 \\
    @VVV @V{\pi}VV @| \\
    SO(6) @>>> SO(7) @>>> S^6
  \end{CD}
\end{equation*}
where the horizontal lines are principal fibre bundles and $p(g) =
\pi(g)e_1$.

Lemma 4.1 of \cite{mn} implies the following lemma immediately.

\begin{lem}
  Put $V^6 = D^3 \times D^2 \times D^1$.
  Then the composite map
  $p \varphi_6 : (V^6, \partial V^6) \rightarrow (S^6, \{e_1\})$
  is a relative homeomorphism.
\end{lem}

Now we can state one of our main results.

\begin{thm}
  \label{thm:spin7}
  The cell complex $e^0 \cup e^3 \cup e^5 \cup e^6 \cup e^7
  \cup e^8 \cup e^9 \cup e^{10} \cup e^{11} \cup e^{12} \cup e^{13}
  \cup e^{14} \cup e^{15} \cup e^{16} \cup e^{18} \cup e^{21}$ gives a
  cellular decomposition of $Spin(7)$.
\end{thm}
\begin{proof}
  First we show that $\ee^i \cap \ee^j
  = \emptyset$ if $i \not= j$.
  We consider the following three cases:

  (1) For the case where $i, j \in \{0, 3, 5, 7, 8, 10, 12, 15 \}$;
  both cells $e^i$ and $e^j$ are in $SU(4)$ and $e^0 \cup e^3 \cup
  e^5 \cup e^7 \cup e^8 \cup e^{10} \cup e^{12} \cup e^{15}$ is a
  cellular decomposition of $SU(4)$,
  whence we have $\ee^i \cap \ee^j
  = \emptyset$ if $i \not= j$.

  (2) For the case where $i \in \{0, 3, 5, 7, 8, 10, 12, 15 \}$ and
  $j \in \{6, 9, 11, 13, 14, 16, 18, \\
  21 \}$;
  we have $p(\ee^i) = \{e_1\}$ and $p(\ee^j)
  = S^6 \backslash \{e_1\}$,
  whence we have $\ee^i \cap \ee^j
  = \emptyset$.

  (3) For the case where $i, j \in \{6, 9, 11, 13, 14, 16, 18, 21 \}$,
  suppose that $A \in \ee^i \cap \ee^j$.
  Since $\ee^i = \ee^6\ee^{i-6}$ and $\ee^j = \ee^6\ee^{j-6}$, we
  can put $A = A_1A_2 = A'_1A'_2$, where $A_1, A'_1 \in \ee^6$, $A_2
  \in \ee^{i-6}$ and $A'_2 \in \ee^{j-6}$.
  We have $A_1 = A'_1$, since $p(A_1) = p(A_1A_2) = p(A'_1A'_2)
  = p(A'_1)$ and $p|_{\ee^6}$ is monic.
  Then we have $A_2 = A'_2$ and the first case shows that $i-6 =
  j-6$, that is, $i=j$.
  Thus $\ee^i \cap \ee^j = \emptyset$ if $i \not=j$.

  Next, we will check that the boundaries of the cells are included in
  the lower dimensional cells.
  In Proposition \ref{prop:su4}, it is proved that
  the boundaries of the cells of $SU(4)$ are included in
  the lower dimensional cells.
  In Proof of Theorem 4.2 in \cite{mn}, we showed that $\dot e^6 \subset
  e^3 \cup e^5$, $\dot e^9 \subset e^6 \cup e^8$, $\dot e^{11} \subset
  e^5 \cup e^9$ and $\dot e^{14} \subset e^8 \cup e^9 \cup e^{11}$.
  By using (2) of Remark \ref{rem:e}, we also obtain
  \begin{align*}
    \dot e^{13} & = e^6\dot e^7 \cup \dot e^6e^7 \subset e^{11} \cup e^{12}, \\
    \dot e^{16} & = e^6e^7\dot e^3 \cup e^6\dot e^7e^3 \cup
    \dot e^6e^7e^3 \subset e^{13} \cup e^{14} \cup e^{15}, \\
    \dot e^{18} & = e^6e^7\dot e^5 \cup e^6\dot e^7e^5 \cup
    \dot e^6e^7e^5 \subset e^{16} \cup e^{14} \cup e^{15}, \\
    \dot e^{21} & = e^6e^7e^5\dot e^3 \cup e^6e^7\dot e^5e^3 \cup
    e^6\dot e^7e^5e^3 \cup \dot e^6e^7e^5e^3 \subset e^{18} \cup
    e^{16} \cup e^{14} \cup e^{15}.
  \end{align*}

  Finally, we will show that the inclusion map $e^0 \cup e^3 \cup e^5
  \cup e^6 \cup e^7 \cup e^8 \cup e^9 \cup e^{10} \cup e^{11} \cup
  e^{12} \cup e^{13} \cup e^{14} \cup e^{15} \cup e^{16} \cup e^{18}
  \cup e^{21} \rightarrow Spin(7)$
  is epic.
  Let $g \in Spin(7)$.
  If $p(g) = e_1$, then $g$ is contained in $SU(4) = e^0 \cup e^3
  \cup e^5 \cup e^7 \cup e^8 \cup e^{10} \cup e^{12} \cup e^{15}$.
  Suppose that $p(g) \not= e_1$.
  There is an element $h\in e^6$ such that $p(h) = p(g)$.
  Thus we have $h^{-1}g \in SU(4)$ since $p(h^{-1}g) = e_1$.
  Therefore we have $g \in h(e^0 \cup e^3 \cup e^5 \cup e^7 \cup e^8
  \cup e^{10} \cup e^{12} \cup e^{15}) \subset
  e^0 \cup e^3 \cup e^5 \cup e^6 \cup e^7 \cup e^8 \cup e^9 \cup
  e^{10} \cup e^{11} \cup e^{12} \cup e^{13} \cup e^{14} \cup e^{15}
  \cup e^{16} \cup e^{18} \cup e^{21}$.
\end{proof}

\begin{rem}
  Araki \cite{araki} also gave a cellular decomposition of $Spin(n)$, 
  but the one we have given here is a cellular decomposition with the
  minimum number of cells, satisfying the Yokota principle
  (\cite{yok1}, \cite{yok2}, \cite{yok3}).
  As will be seen later, it is effectively used to determine the
  Lusternik-Schnirelmann category.
\end{rem}

It is easy to give a cellular decomposition of $Spin(8)$ using a
homeomorphism $Spin(8) \rightarrow Spin(7) \times S^7$.

\section{The cone-decomposition of $SU(4)$}

Obviously there is a filtration $F'_0 = * \subset F'_1 = SU(4)^{(7)} \subset
F'_2 = SU(4)^{(12)} \subset F'_3 = SU(4)$.
It is well-known that $F'_1 = \Sigma \mathbb CP^3 = S^3 \cup e^5
\cup e^7$ and $F'_2 = F'_1 \cup e^8 \cup e^{10} \cup e^{12}$.
Thus the integral cohomology $H^n(F'_2;\mathbb Z)$ is given by
\begin{equation*}
  H^n(F'_2;\mathbb Z) \cong
  \begin{cases}
    \mathbb Z\langle{1}\rangle & (n = 0) \\
    \mathbb Z\langle{y_n}\rangle & (n = 3, 5, 7, 8, 10, 12) \\
    0 & (\text{otherwise}).
  \end{cases}
\end{equation*}
The action of the squaring operation $Sq^2$ is given as follows:
\begin{equation*}
  Sq^2 y_n =
  \begin{cases}
    y_{n+2} & {\rm for} \ n = 3, 10, \\
    0 & {\rm for} \ n = 5, 7, 8, 12
  \end{cases}
\end{equation*}
where $y_n$ is regarded as an element of the mod 2 cohomology.
To give the cone decomposition of $SU(4)$, we use the following
homotopy fibration:
\begin{equation}
  \label{eq:fib}
  F \overset{\Psi}\longrightarrow F'_1 \overset{\iota}\longrightarrow F'_2.
\end{equation}
Without loss of generality, we may regard this as a Hurewicz fibration
over $F'_2$.

Firstly we consider the Serre spectral sequence $(E_r^{*,*}, d_r)$
associated with the above fibration, where the generators of $E_2^{*,0}$
for $* \leq 7$ are permanent cycles and survive to
$E_{\infty}$-terms.
Hence $F$ is $6$-connected and the transgression $\tau : H^7(F;\mathbb
Z) \to H^{8}(F'_2;\mathbb Z)$ is an isomorphism to $H^8(F'_2;\mathbb
Z) \cong {\mathbb Z} \langle {y_{8}}\rangle$.
Thus $H^7(F;\mathbb Z) \cong {\mathbb Z}\langle{x_{7}}\rangle$ for
some $x_7 \in H^7(F;\mathbb Z)$.
Similarly, the generators in $E_2^{3,7}\cong{\mathbb
  Z}\langle{y_{3}{\otimes}x_{7}}\rangle$ and $E_2^{10,0}\cong
H^{10}(F'_2;\mathbb Z) \cong {\mathbb Z}\langle{y_{10}}\rangle$ must
lie in the image of differentials $d_3$ and $d_{10}=\tau :
H^9(F;\mathbb Z) \to H^{10}(F'_2;\mathbb Z)$ respectively, and we have
that $H^8(F;\mathbb Z)=0$ and $H^9(F;\mathbb Z)
\cong {\mathbb Z}\langle{x_{9}}\rangle{\oplus}{\mathbb
  Z}\langle{x'_{9}}\rangle$, where the elements $x_9$ and $x'_9$ in
$H^9(F;\mathbb Z)$ are corresponding to $x_{10}$ and $y_{3}{\otimes}x_{7}$ by
the transgression $\tau$ and $d_3$ respectively.
We remark that the choice of the generator $x'_9$ is not unique.
Continuing this process, we have that $H^{10}(F;\mathbb Z)=0$ and
$H^{11}(F;\mathbb Z)
\cong {\mathbb Z}\langle{x_{11}}\rangle{\oplus}{\mathbb
  Z}\langle{x'_{11}}\rangle{\oplus}{\mathbb
  Z}\langle{x''_{11}}\rangle{\oplus}{\mathbb
  Z}\langle{x'''_{11}}\rangle$ whose generators are corresponding to
$x_{12}$, $y_{3}{\otimes}x_{9}$, $y_{3}{\otimes}x'_{9}$ and
$y_{5}{\otimes}x_{7}$ respectively by the transgression $\tau$ and
differentials $d_3$, $d_3$ and $d_5$.

Thus the integral cohomology $H^n(F;\mathbb Z)$ for $0 \leq n \leq 11$ is given by
\begin{equation*}
  H^n(F;\mathbb Z) \cong
  \begin{cases}
    \mathbb Z\langle{1}\rangle & (n = 0) \\
    \mathbb Z\langle{x_{7}}\rangle & (n = 7) \\
    \mathbb Z\langle{x_{9}}\rangle \oplus \mathbb Z\langle{x'_{9}}\rangle & (n = 9) \\
    \mathbb Z\langle{x_{11}}\rangle \oplus \mathbb Z\langle{x'_{11}}\rangle \oplus \mathbb Z\langle{x''_{11}}\rangle \oplus \mathbb Z\langle{x'''_{11}}\rangle
    & (n = 11) \\
    0 & (\text{otherwise}) 
  \end{cases}
\end{equation*}
where $x_7$, $x_9$ and $x_{11}$ are transgressive generators in
$H^{\ast}(F;\mathbb Z)$.
Hence $F$ has, up to homotopy, a cellular decomposition
$e^0 \cup e^7 \cup_{\varphi_1} e^9 \cup_{\varphi'_1} e^9_1
\cup_{\varphi_2} e^{11} \cup$ (cells in dimensions $\geq 11$), where the
cells $e^7$, $e^9$ and $e^{11}$ correspond to $x_7$, $x_9$ and
$x_{11}$ respectively.
Then we obtain a subcomplex $A' = e^0 \cup e^7 \cup_{\varphi_1} e^9
\cup_{\varphi'_1} e^9_1 \cup_{\varphi_2} e^{11}$ of $F$.

Secondly, we determine the attaching maps $\varphi_1$ and $\varphi'_1$:
Let us recall that $\pi_8(S^7) \cong {\mathbb Z/2} \langle \eta_7
\rangle$ whose generator $\eta_7$ can be detected by $Sq^2$, the
mod $2$ Steenrod operation.
Since the action of mod $2$ Steenrod operation commutes with the
cohomology transgression (see \cite[Proposition 6.5]{mc}),
we see that $Sq^2x_7$ is transgressive, and
hence is $cx_9$ for some $c \in {\mathbb Z/2}$.
We know that $\tau{x_9} = y_{10} \not= 0$ and $\tau Sq^2x_7 =
Sq^2\tau{x_7} = Sq^2y_8 = 0$, and hence $Sq^2x_7$ must be trivial.
Thus the attaching maps $\varphi_1$ and $\varphi'_1$ are both null
homotopic and $A'$ is homotopy equivalent to $(S^7 \vee S^9 \vee S^9_1)
\cup_{\varphi_2} e^{11}$.

Thirdly we check the composition of projections with the attaching map
$\varphi_2:S^{10} \rightarrow S^7 \vee S^9 \vee S^9_1$ to $S^9$ and
$S^9_1$, which can also be detected by $Sq^2$.
Again by the commutativity of the action of mod 2 Steenrod
operation with the transgression, we see that the composition map
$\proj_{S^9} \circ \varphi_2 : S^{10}
\overset{\varphi_2}{\longrightarrow} S^7 \vee S^9  \vee S^9_1
\longrightarrow S^9$
represents a generator of $\pi_{10}(S^9) \cong \mathbb
Z/2\langle\eta_{9}\rangle$, since $Sq^2 : H^8(F'_2;\mathbb Z/2)
\rightarrow H^{10}(F'_2;\mathbb Z/2)$ is non-trivial.
If the composition map
$\phi_1 = \proj_{S^9_1} \circ \varphi_2 : S^{10} \overset{\varphi_2}
{\longrightarrow} S^7 \vee S^9  \vee S^9_1 \longrightarrow S^9_1$
is non-trivial, we replace $\varphi_2$ by the composition of
$\varphi_2$ and the homotopy equivalence $\xi : S^7 \vee S^9 \vee
S^9_1 \to S^7 \vee S^9 \vee S^9_1$ where $\xi|_{S^7}$ and
$\xi|_{S^9_1}$ are the identity maps and $\xi|_{S^9}$ is the unique
co-H-structure map $\phi : S^9 \rightarrow S^9 \vee S^9_1$;
then we obtain that $\phi_1$ is trivial, since $2\eta_9 = 0$.
Then $A'$ is homotopy equivalent to $((S^7 \vee S^9) \cup_{\varphi_2}
e^{11}) \vee S^9_1$.
Let $A$ denote the subcomplex $(S^7 \vee S^9) \cup_{\varphi_2} e^{11}$ of $A'$ and $\psi = \Psi\vert_{A} : A \to F'_1$.

\begin{lem}
  \label{lem:3.1}
  $F'_2$ is homotopy equivalent to $F'_1 \cup_{\psi} CA$.
\end{lem}
\begin{proof}
  The image of $H^*(A;{\mathbb Z})$ in $H^*(F;{\mathbb Z})$ under the
  induced map of the inclusion coincides with the module of
  transgressive elements with respect to the fibration (\ref{eq:fib})
  (see \cite[Chapter 6]{mc}).
  Thus we may regard that 
  $H^{n-1}(A;{\mathbb Z}) = \delta^{-1}(\iota^{\ast}(H^n(F'_2,\ast;))) \subset 
  H^{n-1}(F;{\mathbb Z})$:
  \begin{equation*}
    \begin{CD}
      H^{n-1}(F;\mathbb Z) @>{\delta_F}>> 
      H^n(F'_1,F;{\mathbb Z}) @<{\iota_F^*}<< 
      H^n(F'_2,\ast; {\mathbb Z})
      \\
      @V{}VV @V{}VV @V{\Vert}VV 
      \\
      H^{n-1}(A;\mathbb Z) @>{\delta_A}>> 
      H^n(F'_1,A;{\mathbb Z}) @<{\iota_A^*}<< 
      H^n(F'_2,\ast; {\mathbb Z}),
    \end{CD}
  \end{equation*}
  where $\iota_F$ and $\iota_A$ are given by $\iota$, and $\delta_F$ and 
  $\delta_A$ denote the connecting homomorphisms of the long exact sequences 
  for the pairs $(F'_1,F)$ and $(F'_1,A)$, respectively.
  Thus the image of $\delta_A$ is contained in the image of $\iota_A^{\ast}$ 
  and we also have 
  \begin{equation*}
    H^n(F'_1,A;{\mathbb Z}) \cong H^n(F'_1 \cup_{\psi} CA,CA;{\mathbb Z})
    \cong H^n(F'_1 \cup_{\psi} CA,\ast;{\mathbb Z}).
  \end{equation*}
  Since the composition map $A \overset{\psi}\rightarrow F'_1
  \overset{\iota}\rightarrow F'_2$ is trivial, we can define a map
  \begin{equation*}
    f : F'_1 \cup_{\psi} CA \longrightarrow F'_2,
  \end{equation*}
  by $f|_{F'_1} = \iota : F'_1 \rightarrow F'_2$ and $f|_{CA} = \ast$.

  In order to prove the lemma, we show that 
  $f^* : H^n(F'_2;{\mathbb Z}) \cong {\mathbb Z} \rightarrow 
  H^n(F'_1\cup_{\psi} CA;{\mathbb Z}) \cong {\mathbb Z}$
  is an isomorphism for $n = 3, 5, 7, 8, 10, 12$.
  We have a commutative diagram
  \begin{equation*}
    \begin{CD}
      @. H^n(F'_2;\mathbb Z) @>{\iota^*}>> H^n(F_1';\mathbb Z) \\
      @. @V{f^*}VV @V{\Vert}VV \\
      H^n(F'_1 \cup CA, F'_1;\mathbb Z) @>{j^*}>>
      H^n(F'_1 \cup CA;\mathbb Z) @>{i^*}>> H^n(F'_1;\mathbb Z),
    \end{CD}
  \end{equation*}
  where the bottom row is a part of the exact sequence for the
  pair $(F'_1 \cup CA, F'_1)$.
  The induced map $i^*$ is an isomorphism for $n \leq 7$, since
  $H^n(F'_1 \cup CA, F'_1;\mathbb Z) = 0$ for $n \leq 7$ and
  since $\iota^*$ is an isomorphism for $n \leq 7$.
  Then we obtain that $f^*$ is an isomorphism for $n \leq 7$.
  Moreover we can show that $j^* : H^n(F'_1 \cup CA, F'_1;\mathbb Z)
  \rightarrow H^n(F'_1 \cup CA;\mathbb Z)$ is an isomorphism for $n
  \geq 8$, by considering the exact sequence for the pair $(F'_1 \cup
  CA, F'_1)$, since we have $H^n(F'_1)=0$ for $n \geq 8$.
  To perform the other cases for $n = 8, 10, 12$, it is sufficient to
  show that $f^*$ is surjective.
  In fact, we have a commutative diagram

  \xymatrix{
    & H^{n-1}(A;\mathbb Z) \ar[r]^{\delta_A} \ar[d]^{\cong}_{\Sigma}
    & H^n(F'_1,A;\mathbb Z) \ar[rd]^{\cong}
    & H^n(F'_2,*;\mathbb Z) \ar[l]_-{\iota_A^*} \ar[d]^{f^*} \\
    & H^n(\Sigma A,*;\mathbb Z) \ar[r]^-{\cong}
    & H^n(F'_1 \cup CA, F'_1;\mathbb Z) \ar[u] \ar[r]^-{j^*}
    & H^n(F'_1 \cup CA, *;\mathbb Z),
    }

  \medskip
  \noindent
  where $\Sigma$ is the suspension isomorphism.
  Since $j^*$ is an isomorphism for $n \geq 8$, we obtain that
  $\delta_A$ is an isomorphism for $n \geq 8$.
  Since the image of $\delta_A$ is contained in the image $\iota_A^*$,
  we see that $f^*$ is surjective for $n \geq 8$, and hence $f$ is a
  homotopy equivalence.
\end{proof}

\begin{prop}
  We have $w{\rm cat}(F'_i) = {\rm cat}(F'_i) =
  {\rm Cat}(F'_i) = i$.
\end{prop}
\begin{proof}
  The cohomology of $F'_i$ implies that $w{\rm cat}(F'_i) \geq i$.
  The cone-decomposition
  \begin{equation*}
    F'_1 = \Sigma \mathbb CP^3, \quad F'_2 \simeq F'_1 \cup CA,
    \quad F'_3 = F'_2 \cup CS^{14}
  \end{equation*}
  implies that ${\rm Cat}(F'_i) \leq i$, which completes the proof.
\end{proof}

\section{Proof of Theorem \ref{thm:main}}

We define a filtration $F_0 = * \subset F_1 \subset F_2 \subset F_3
\subset F_4 \subset F_5 =Spin(7)$ by
\begin{equation*}
  \begin{array}{ll}
    F_1 = SU(4)^{(7)}, &
    F_2 = SU(4)^{(12)} \cup e^6, \\
    F_3 = SU(4) \cup e^6 \cup e^9 \cup e^{11} \cup e^{13}, &
    F_4 = Spin(7)^{(18)}.
  \end{array}
\end{equation*}

We need the following lemma to prove Theorem \ref{thm:filtration}.
\begin{lem}
  \label{lem:join}
  We have a homeomorphism of pairs
  \begin{equation*}
    (CA_1, A_1) \times (CA_2, A_2) = (C(A_1 * A_2), A_1 * A_2).
  \end{equation*}
\end{lem}
(The proof can be found in p.482-483 of \cite{white2}.)

Now Theorem \ref{thm:main} follows from the following theorem.

\begin{thm}
  \label{thm:filtration}
  We have $w{\rm cat}(F_i) = {\rm cat}(F_i) =
  {\rm Cat}(F_i) = i$.
\end{thm}
\begin{proof}
  The mod 2 cohomology of $F_i$ implies that $w{\rm cat}(F_i) \geq i$.
  Then it is sufficient to show that ${\rm Cat}(F_i) \leq i$.
Obviously we have a homeomorphism $F_1 = \Sigma \mathbb CP^3$.
Since the cell $e^6$ is attached to $F_1$,
we obtain that $F_2 \simeq F_1 \cup C(S^5 \vee A)$ using Lemma
\ref{lem:3.1}.
Since we have $e^9 \cup e^{11} \cup e^{13} = e^6(e^3 \cup e^5 \cup
e^7)$, the composition map
\begin{equation*}
  (CS^5, S^5) \times (C\mathbb CP^3, \mathbb CP^3) \longrightarrow
  (CS^5, S^5) \times (\Sigma \mathbb CP^3, *) \longrightarrow
  (F_2 \cup e^9 \cup e^{11} \cup e^{13}, F_2)
\end{equation*}
is a relative homeomorphism.
Then we obtain $F_2 \cup e^9 \cup e^{11} \cup e^{13} = F_2 \cup C(S^5
* \mathbb CP^3)$ using Lemma \ref{lem:join}.
The cell $e^{15}$ is the highest dimensional cell of $SU(4)$ and is
attached to $F_2$.
Then we obtain $F_3 \simeq F_2 \cup C(S^{14} \vee (S^5 * \mathbb
CP^3))$.
Now we consider the following composition map:
\begin{equation*}
  (C(S^5*A), S^5*A) = (CS^5, S^5) \times (CA, A) \longrightarrow
  (CS^5, S^5) \times (F'_2, F'_1) \longrightarrow (F_4, F_3).
\end{equation*}
Since we have $e^{14} \cup e^{16} \cup e^{18} = e^6 (e^8 \cup e^{10}
\cup e^{12})$, the right map is a relative homeomorphism.
The left map induces an isomorphism of homologies of pairs
so that the map $H_*(F_3 \cup C(S^5*A), F_3 ; \mathbb Z)
\rightarrow H_*(F_4, F_3 ; \mathbb Z)$ is an isomorphism.
Thus we obtain $F_4 \simeq F_3 \cup C(S^5*A)$.
Obviously we have a homeomorphism $F_5 = F_4 \cup CS^{20}$.
\end{proof}

\end{document}